\documentclass{conm-p-l}

\copyrightinfo{2006}{}

\newtheorem{theorem}{Theorem}[section]
\newtheorem{lemma}[theorem]{Lemma}

\theoremstyle{definition}
\newtheorem{definition}[theorem]{Definition}

\theoremstyle{remark}

\numberwithin{equation}{section}

\newcommand{\tu}{\tilde{u}}

\newcommand{\cS}{{\mathcal S}}

\newcommand{\e}{\epsilon}

\newcommand{\A}{{\mathcal A}}

\newcommand{\Ga}{\Gamma}

\newcommand{\dl}{\delta}

\newcommand{\Sg}{\Sigma}

\newcommand{\z}{\zeta}

\newcommand{\na}{\nabla}
\newcommand{\lag}{\langle}
\newcommand{\rag}{\rangle}

\newcommand{\vphi}{\varphi}

\newcommand{\tB}{\tilde{B}}

\newcommand{\ta}{\tilde{a}}
\newcommand{\tn}{\tilde{n}}
\newcommand{\ha}{\hat{a}}
\newcommand{\hn}{\hat{n}}
\newcommand{\tp}{\tilde{p}}

\def\maprightu#1{\smash{
    \mathop{\longrightarrow}\limits^{#1}}}
\def\maprightd#1{\smash{
    \mathop{\longrightarrow}\limits_{#1}}}
\def\mapdownl#1{
    \llap{$\vcenter{\hbox{$\scriptstyle#1$}}$}\Big\downarrow}
\def\mapdownr#1{\Big\downarrow
    \rlap{$\vcenter{\hbox{$\scriptstyle#1$}}$}}

\begin{document}

\title[Segment Description of Turbulence]
{Segment Description of Turbulence}

\author{Y. Charles Li}
\address{Department of Mathematics, University of Missouri, 
Columbia, MO 65211, USA}
\curraddr{}
\email{cli@math.missouri.edu}
\thanks{}


\subjclass{Primary 35, 37; Secondary 34, 46}
\date{}

\dedicatory{}


\begin{abstract}
We propose a segment description for turbulent solutions to evolution equations 
in an effort to develop an effective description rather than the classical Reynolds 
average. The new description has connections with symbolic dynamics and shadowing 
technique. The challenge of future study is how to effectively implement the description 
numerically. 
\end{abstract}

\maketitle










\section{Introduction}

Chaos and turbulence have no good averages \cite{Li06i}. The matter is more fundamental 
than just poor understanding of averages. The very mechanism of chaos leads to the 
impossibility of a good average \cite{Li04}. From what we learn about chaos in partial 
differential equations \cite{Li04}, turbulent 
solutions not only have sensitive dependences on initial conditions, but 
also are densely packed inside a domain in the phase space. They are far 
away from the feature of fluctuations around a mean. In fact, they wander 
around in a fat domain rather than a thin domain in the phase space. 
Therefore, averaging makes no sense at all. One has to seek other descriptions. 
In this article, we present a description -- named segment description. 
Segment description approximates a turbulent solution by an infinite sequence of segments
over an infinite time interval. The total number of different segments is finite, and all
the segments are over a fixed finite time interval. Thus, each segment can be calculated 
rather accurately. The segment description also has connections with both symbolic 
dynamics and shadowing technique \cite{Li04}. But, our emphasis here is upon applications 
rather than mathematical abstract arguments. In our view, there are two most significant 
problems in the study of turbulence: 1. an effective description of turbulent solutions 
rather than the Reynolds average, 2. control of turbulence. Here the segment description 
is aiming at problem 1 by offering an alternative description. 

In a sense, the Reynolds average is based upon the feature of fluctuations around a mean, 
which turbulent solutions do not have. Consider the Navier-Stokes equations in the form:
\begin{equation}
u_{i,t} + u_j u_{i,j} = - p_{,i} + \ \mbox{Re}^{-1} \ u_{i,jj} +f_i \ , 
\quad u_{i,i} = 0 \ ;
\label{NS}
\end{equation}
where $u_i$'s are the velocity components, $p$ is the pressure, $f_i$'s 
are the external force components, and Re is the Reynolds number.
According to the Reynolds average, one splits the 
velocity and pressure of fluids into two parts:
\[
u_i = U_i + \tu_i \ , \quad p=P+\tp
\]
where the capital letters represent relatively long time averages which 
are still a function of time and space, and the tilde-variables represent 
mean zero fluctuations,
\[
U_i = \lag u_i \rag \ , \quad \lag \tu_i \rag =0 \ , \quad
P = \lag p \rag \ , \quad \lag \tp \rag =0 \ .
\]
A better interpretation is by using ensemble average of repeated 
experiments. One can derive the Reynolds equations for the averages,
\begin{equation}
U_{i,t} + U_jU_{i,j} = - P_{,i} + \ \mbox{Re}^{-1} \ U_{i,jj} 
-\lag \tu_i\tu_j \rag_{,j} +f_i \ , 
\quad U_{i,i} = 0 \ .
\label{RE}
\end{equation}
The term $\lag \tu_i\tu_j \rag$ is completely unknown. Fluid engineers 
call it Reynolds stress. The Reynolds model is given by
\begin{equation}
\lag \tu_i\tu_j \rag = -R^{-1} \ U_{i,j} \ ,
\label{RM}
\end{equation}
where $R$ is a constant. There are many more models on the term 
$\lag \tu_i\tu_j \rag$ \cite{Hin75}. But no one leads to a satisfactory result. One 
can re-interpret the Reynolds equations (\ref{RE}) as control equations 
of the orginal Navier-Stokes equations (\ref{NS}), with the term 
$\lag \tu_i\tu_j \rag_{,j}$ being the control of taming a turbulent solution 
to a laminar solution (hopefully nearby) \cite{Li06i}. The Reynolds model (\ref{RM}) 
amounts to changing the fluid viscosity which can bring a turbulent flow 
to a laminar flow. This laminar flow may not be anywhere near the turbulent 
flow though. Thus, the Reynolds model may not produce satisfactory result 
in comparison with the experiments. Fluid engineers gradually gave up all 
these Reynolds' type models and started directly computing the original 
Navier-Stokes equations (\ref{NS}).

\section{Segment Description}

Let $\cS$ be the phase space which is a metric space, on which a flow $F^t$ ($t \geq 0$) is 
defined. $F^t$ is $C^0$ in $t$, and for any fixed $t \in (0, +\infty )$, $F^t$ is a $C^1$ 
map. For instance, $\cS$ can be a Sobolev space defined on the 2-tori, and $F^t$ is the 
2D Navier-Stokes flow. Let $\A$ be a compact absorbing set on which the dynamics of $F^t$ 
can be turbulent. The attractor is a subset of $\A$. $\A$ can be chosen to be just the 
attractor. But a compact absorbing set is easier to obtain than an attractor in numerics.
Let $B$ be a ball that includes $\A$, such that
\begin{equation}
\| DF^t(x) \| \leq C(t), \quad \forall x \in B, 
\label{derb}
\end{equation}
where $C(t)$ is a continuous function of $t \in [0, +\infty )$. 
Let $D$ be any subset of $\cS$, the diameter of $D$ is defined as
\[
\text{diameter}\{ D\} = \sup_{x,y \in D} \text{distance}(x,y)\ .
\]
Let $B_\dl (x)$ be the ball centered at $x$ with radius $\dl >0$. Fix a $T >0$ and a
small $\e >0$, for any $x \in \A$, there exists a $\dl >0$ [e.g. $\dl = c \e$ for some 
$c>0$ by (\ref{derb})] such that 
\begin{equation}
\sup_{t\in [0,T]} \text{diameter}\left \{ F^t (B_\dl (x))\right \} \leq \e \ .
\label{dir}
\end{equation}
The union of all such balls 
\[
\bigcup_{x \in \A } B_\dl (x)
\]
is a covering of $\A$. Since $\A$ is compact, there is a finite covering
\[
\bigcup_{1\leq n \leq N} B_{\dl_n} (x_n) \ .
\]
We can choose this covering to be minimal. From this covering, we can deduce a mutually 
disjoint covering:
\begin{eqnarray*}
\tB_{\dl_n} (x_n) &=& B_{\dl_n} (x_n) \bigcap 
\overline{\bigcup_{n+1\leq m \leq N} B_{\dl_m} (x_m)} \ , \quad (n=1, \cdots , N-1)\ ; \\
\tB_{\dl_N} (x_N) &=& B_{\dl_N} (x_N)\ .
\end{eqnarray*}
Then $\bigcup_{1\leq n \leq N} \tB_{\dl_n} (x_n)$ is a mutually disjoint covering of $\A$.
Thus 
\[
\left \{ A_n = \tB_{\dl_n} (x_n) \cap \A \ , \quad (1\leq n \leq N) \right \} 
\]
form a partition of $\A$. Let $\mu$ be a Borel measure defined on $\A$, then the metric 
entropy of the partition is 
\[
H = -\sum_{n=1}^N \mu (A_n) \log \mu (A_n) \ .
\]
In this article, we are not going to emphasize either metric entropy or topological entropy
\cite{You03} since they are essentially sums of positive Liapunov exponents and Liapunov 
exponents are very easy to compute in numerics. Unlike the usual partitions, here $A_n$'s 
are of small size. On the other hand, we do not take the size to the zero limit \cite{You03}
in which case we will be led to Liapunov exponents. Our goal is to obtain a uniform 
approximation of the turbulent solutions over an infinite time interval.
\begin{definition}[Segments]
Define the segments $s_n$ as follows
\[
s_n = \bigcup_{t\in [0,T]} F^t(x_n)\ , \quad (n=1,\cdots , N)\ . 
\]
\end{definition}
\begin{definition}[Segment Description]
For any $x \in \A$, the orbit
\begin{equation}
\xi (t) = \left \{ F^t(x) \ , \quad t\in [0, +\infty ) \right \} 
\label{orb}
\end{equation}
has the segment representation
\begin{equation}
\eta (t) = \left \{ s_{n_0} s_{n_1} s_{n_2} \cdots \right \} \ , \quad t\in [0, +\infty )
\label{sgr}
\end{equation}
where 
\[
F^{jT}(x) \in A_{n_j}\ , \quad n_j \in \{ 1,\cdots , N \} \ , \quad 0\leq j <+\infty \ .
\]
We call $\eta (t)$ an $\e$-pseudo-orbit.
\label{sr}
\end{definition}
\begin{lemma}
For any $t\in [0, +\infty )$, 
\[
\text{distance} \left ( \xi (t), \eta (t) \right )  \leq \e \ .
\]
\label{edi}
\end{lemma}
\begin{proof}
The lemma follows from the condition (\ref{dir}) and the Definition \ref{sr}.
\end{proof}
We define a {\em maximal difference function} as follows
\[
M_d(t) = \max_{1 \leq m,n \leq N} \text{distance} \left ( s_m (t), s_n (t) \right )\ .
\]
If the absorbing set $\A$ is actually an attractor, then $M_d(t)$ should be almost 
constant in $t$ and almost the diameter of the attractor. In general, $M_d(t)$ 
serves as a measure of the shrinking of the diameter of the absorbing set $\A$ 
under the flow.

\subsection{Markov Description \label{mard}}

For any $1\leq m,n \leq N$, if 
\begin{equation}
F^T(A_m) \cap A_n \neq \emptyset \ ,
\label{inc}
\end{equation}
then we assign the value $\Ga_{mn}=1$; otherwise $\Ga_{mn}=0$. These $\Ga_{mn}$'s form 
a $N\times N$ matrix $\Ga$ called the transition matrix. In the segment representation 
(\ref{sgr}), the first segment $s_{n_0}$ can be any of the $N$ choices, while the second 
segment $s_{n_1}$ is admissible only when $\Ga_{n_0n_1}=1$. The admissibility of the 
third segment $s_{n_2}$ depends on both $s_{n_0}$ and $s_{n_1}$. In general, the 
admissibility of the $j$-th segment $s_{n_{j-1}}$ depends on all segments in front of it
$s_{n_0} \cdots s_{n_{j-2}}$. Thus we are dealing with a non-Markov process. Nevertheless, 
we believe that the Markov description is still significant. The Markov description extended
the admissibility by the criterion that for any $j=1,2, \cdots$; $s_{n_j}$ is admissible if 
$\Ga_{n_{j-1}n_j}=1$. We can define a probability as follows
\begin{equation}
p_{mn} = \mu \left ( F^T(A_m) \cap A_n \right ) \bigg / \mu (F^T(A_m)) \ ,
\label{mm}
\end{equation}
where $\mu$ is the Borel measure defined on $\A$. 
\begin{definition}
Let $\Xi_1$ be the space consisting of elements of the form
\[
a = \{ n_0 n_1 n_2 \cdots \}
\]
where $n_0 \in \{ 1, \cdots , N \}$ and the admissibility of $n_j$ is decided by 
$\Ga_{n_{j-1}n_j}=1$ for $j=1,2,\cdots$. The product topology in $\Xi_1$ is given by taking 
as the neighborhood basis of any element
\[
a^* = \{ n^*_0 n^*_1 n^*_2 \cdots \}
\]
the cylinder sets
\[
W_j = \left \{ a \in \Xi_1 \ | \ n_k = n^*_k \quad (k <j) \right \}
\]
for $j=1,2,\cdots$. The left shift $\chi_1$ is defined by 
\[
\chi_1(a) = \{ n_1 n_2 n_3 \cdots \} \ , \quad \forall a = \{ n_0 n_1 n_2 \cdots \} \ .
\]
The pair ($\Xi_1, \chi_1$) is called a subshift of finite type.
\label{symd}
\end{definition}
\begin{lemma}
If each row of $\Ga$ has more than one entry to be $1$, then the left shift $\chi_1$ has 
sensitive dependence on the initial conditions.
\label{sle}
\end{lemma}
\begin{proof}
For any large $j$, let $\ta , \ha \in W_j$ and 
\[
\ta = \{ \tn_0 \tn_1 \tn_2 \cdots \} \ , \quad 
\ha = \{ \hn_0 \hn_1 \hn_2 \cdots \} \ ;
\]
then 
\[
\tn_k = \hn_k \quad (k<j) \ .
\]
Since each row of $\Ga$ has more than one entry to be $1$, we can choose such $\ta$ and $\ha$
that $\tn_j \neq \hn_j$. After $j$ iterations of $\chi_1$, we have
\[
\chi^j_1(\ta )= \{ \tn_j \tn_{j+1} \tn_{j+2} \cdots \} \ , \quad 
\chi^j_1(\ha )= \{ \hn_j \hn_{j+1} \hn_{j+2} \cdots \} \ .
\]
Thus
\[
\chi^j_1(\ha ) \not\in W_1 \text{ of } \chi^j_1(\ta ) \ .
\]
That is, after $j$ iterations of $\chi_1$, the two elements are far away. The proof is 
complete. 
\end{proof}
In the case that not all the rows of $\Ga$ have more than one entry to be $1$, the 
left shift $\chi_1$ may exhibit transient chaos (intermittence) or regular dynamics.
\begin{definition}
Let $\Xi_2$ be the space consisting of elements of the form
\[
\z (t) = \{ s_{n_0} s_{n_1} s_{n_2} \cdots \}
\]
where $n_0 \in \{ 1, \cdots , N \}$ and the admissibility of $n_j$ is decided by 
$\Ga_{n_{j-1}n_j}=1$ for $j=1,2,\cdots$. The product topology in $\Xi_1$ induces a 
product topology in $\Xi_2$ with the neighborhood basis of any element
\[
\z^* (t) = \{ s_{n^*_0} s_{n^*_1} s_{n^*_2} \cdots \}
\]
the cylinder sets
\[
V_j = \left \{ \z (t) \in \Xi_2 \ | \ n_k = n^*_k \quad (k <j) \right \}
\]
for $j=1,2,\cdots$. The left shift $\chi_2$ is defined by 
\[
\chi_2(\z (t)) =  \{ s_{n_1} s_{n_2} s_{n_3} \cdots \} \ , \quad \forall 
\z (t) = \{ s_{n_0} s_{n_1} s_{n_2} \cdots \} \ .
\]
\label{segn}
\end{definition}
\begin{definition}
Let $\Sg$ be the space consisting of orbits starting from any point in $\A$
\[
\xi (t) = \left \{ F^t(x) \ , \quad t\in [0, +\infty ) \right \} \ , \quad 
\forall x \in \A \ .
\]
The product topology in $\Xi_2$ induces a topology in $\Sg$ via the map 
$\vphi \ : \ \Sg \mapsto \Xi_2$ which maps (\ref{orb}) to (\ref{sgr}). For any 
$\xi^* (t) \in \Sg$, $\xi^* (t)$ has a segment representation as given in 
(\ref{sgr})
\[
\eta^* (t) = \left \{ s_{n^*_0} s_{n^*_1} s_{n^*_2} \cdots \right \} \ , 
\quad t\in [0, +\infty )
\]
at which the neighborhood basis $V_j$'s are defined in Definition \ref{segn}. 
The preimages of the $V_j$'s under $\vphi$ form the induced neighborhood basis
at $\xi^* (t)$
\[
U_j = \vphi^{-1}(V_j) \quad \quad (j=1,2,\cdots )\ .
\]
\label{orn}
\end{definition}
In general, the map $\vphi \ : \ \Sg \mapsto \Xi_2$ which maps (\ref{orb}) to (\ref{sgr})
is neither 1-1 nor onto. Under the topologies defined in Definitions \ref{segn} and
\ref{orn}, $\vphi$ is a continuous map.
\begin{theorem}
There is a continuous map $\phi \ : \ \Sg \mapsto \Xi_1$ such that the following diagram 
commutes:
\[
\begin{array}{ccc}
\Sg &\maprightu{\phi} & \Xi_1 \\
\mapdownl{F^T} & & \mapdownr{\chi_1}\\
\Sg & \maprightd{\phi} & \Xi_1
\end{array} 
\]
\label{md}
\end{theorem}
\begin{proof}
Let $\psi$ be the homeomorphism 
\[
\psi \ : \ \Xi_2 \mapsto \Xi_1 \ , \quad \psi (\{ s_{n_0} s_{n_1} s_{n_2} \cdots \}) = 
\{ n_0 n_1 n_2 \cdots \} \ .
\]
Then 
\[
\phi = \psi \circ \vphi \ : \ \Sg \mapsto \Xi_1 
\]
is a continuous map. For any $\xi (t) \in \Sg$, $\xi (t)$ has a segment representation given by
(\ref{sgr}),
\[
\eta (t) = \left \{ s_{n_0} s_{n_1} s_{n_2} \cdots \right \} \ .
\]
Then $F^T(\xi (t))$ has the segment representation
\[
\left \{ s_{n_1} s_{n_2} s_{n_3} \cdots \right \} \ .
\]
Thus the diagram commutes.
\end{proof}
Theorem \ref{md} serves as a Markov description of ($\Sg , F^T$). Notice that the probability
(\ref{mm}) generates a Markov measure on the subshift of finite type ($\Xi_1,\chi_1$). The 
Markov measure of the cylinder sets $W_j$ in Definition \ref{symd} is given by
\[
\mu (W_j) = p_{n^*_0n^*_1}p_{n^*_1n^*_2} \cdots p_{n^*_{j-2}n^*_{j-1}}\ .
\]
The Kolmogorov-Sinai entropy is 
\[
H_\mu = -\sum_{m=1}^N \sum_{n=1}^N p_{mn} \log p_{mn} \ .
\]
Using the Markov description, one can estimate interesting quantities. Let $Q(x)$ be an 
interesting quantity (for example dissipation rate or energy for Navier-Stokes equations), 
one can calculate
\[
q^s_n = \sup_{t\in [0,T]} Q(s_n(t))\ , \quad q^i_n = \inf_{t\in [0,T]} Q(s_n(t))\ ,
\quad  \quad  (n=1,\cdots ,N)\ .
\]
For any finite time interval [$0, mT$] ($m \in \mathbb{Z}^+$), let $x_0 \in A_{n_0}$, then 
one can track all the possible finite sequences starting from $n_0$
\begin{equation}
\{ n_0n_1n_2 \cdots n_{m-1}\} 
\label{pfs}
\end{equation}
by the condition $\Ga_{n_{j-1}n_j}=1$ ($j=1,\cdots , m-1$) (cf. Definition \ref{symd}). In 
particular, one can identify the subset $I^m_{n_0}$ of $\{ 1,\cdots ,N \}$ that consists of 
elements appearing in the sequences (\ref{pfs}). Then one can estimate $Q(\xi (t))$ ($\xi (0)
=x_0$) over the time interval [$0, mT$] by 
\[
q^s = \max_{n \in I^m_{n_0}} q^s_n \ , \quad q^i = \max_{n \in I^m_{n_0}} q^i_n \ .
\]
That is
\[
q^i \leq Q(\xi (t)) \leq q^s \ , \quad t \in [0, mT]\ .
\]
One should expect that numerically it should be quite easy to track the sequences (\ref{pfs}) 
up to very large $m$.

\subsection{Non-Markov Description}

In general, not all the elements in $\Xi_2$ can be realized by an orbit in $\Sg$. Thus all 
the orbits in $\Sg$ are labeled by a subset of $\Xi_2$ via (\ref{sgr}). To better describe the 
segment description (\ref{orb})-(\ref{sgr}), we introduce the transition tensors
\[
\Ga_{n_0n_1} \ , \quad \Ga_{n_0n_1n_2} \ , \quad \cdots 
\]
where $\Ga_{n_0n_1}$ is defined by (\ref{inc}), if 
\[
F^T(F^T(A_{n_0}) \cap A_{n_1})\bigcap A_{n_2} \neq \emptyset \ ,
\]
then $\Ga_{n_0n_1n_2}=1$, otherwise $\Ga_{n_0n_1n_2}=0$. Similarly one can define the rest 
of the transition tensors.
\begin{definition}
We say that the transition tensors are expanding if for any $\Ga_{n_0n_1 \cdots n_j}=1$, there 
is an integer $m>0$ such that there are at least two different values of $n_{j+m}$ for which 
$\Ga_{n_0n_1 \cdots n_j \cdots n_{j+m}}=1$.
\end{definition}
\begin{theorem}
If the transition tensors are expanding, then $F^T \ : \ \Sg \mapsto \Sg$ has sensitive 
dependence on initial conditions. That is, $F^T$ is chaotic.
\label{bct}
\end{theorem}
\begin{proof}
For any $\xi (t) \in \Sg$ and any neighborhood $U_j$ (Definition \ref{orn}), since 
the transition tensors are expanding, there is another element $\xi^* (t) \in U_j$ 
such that the segment representations of $\xi (t)$ and $\xi^* (t)$ are of the form
\[
\left \{ s_{n_0} s_{n_1} s_{n_2} \cdots s_{n_j} \cdots s_{n_{j+m}} \cdots \right \} 
\text{ and } 
\left \{ s_{n_0} s_{n_1} s_{n_2} \cdots s_{n_j} \cdots s_{n^*_{j+m}} \cdots \right \} 
\]
where $n_{j+m} \neq n^*_{j+m}$. Then 
\[
F^{(j+m)T}(\xi^* (t)) \not\in U_1 \text{ of } F^{(j+m)T}(\xi (t)) \ .
\]
Thus $F^T$ has sensitive dependence on initial conditions. 
\end{proof}
If the transition tensors are not expanding, then $F^T$ may exhibit transient chaos 
(intermittence) or regular dynamics. 
\begin{definition}
Let $\Xi_*$ be the space consisting of elements of the form
\[
a = \{ n_0 n_1 n_2 \cdots \}
\]
where for any $j =0,1,2,\cdots$, $n_j \in \{ 1, \cdots , N \}$; and $\Ga_{n_0 \cdots 
n_k}=1$ for any $k=1,2,\cdots$.
The product topology in $\Xi_*$ and the left shift $\chi$ are defined as in 
Definition \ref{symd}.
\end{definition}
It is easy to see that $\Xi_*$ is the image of $\phi$ in Theorem \ref{md}. Thus as a 
corollary of Theorem \ref{md}, we have the following theorem.
\begin{theorem}
There is a continuous map $\phi \ : \ \Sg \mapsto \Xi_*$ which is onto, such that 
the following diagram commutes:
\[
\begin{array}{ccc}
\Sg &\maprightu{\phi} & \Xi_* \\
\mapdownl{F^T} & & \mapdownr{\chi}\\
\Sg & \maprightd{\phi} & \Xi_*
\end{array} 
\]
\end{theorem}

\section{Numerical Implementation}

In a numerical implementation, one starts with a compact absorbing set $\A$ (for 
example, a ball in a subspace spanned by a finite number of Fourier modes) in a 
Banach space. Then one can collocate as many points as possible in the absorbing set $\A$.
Starting from these points, one can generate the segments over a time interval 
\[
s_n(t)\ , \quad (n=1,\cdots ,N)\ , \quad t \in [0,T]\ .
\]
One can also compute the maximal difference function
\[
M_d(t)=\max_{1 \leq m,n\leq N} \| s_m(t)-s_n(t)\| \ .
\]
To implement a Markov description, one can compute the norm of the gradient of $F^T$ 
restricted to $\A$ at $s_n(T)$:
\[
\rho_n = \left \| \na F^T (s_n(T))|_\A \right \| \ .
\]
Let 
\[
r_n = \min_{1 \leq j\leq N, j\neq n} \{ \| s_n(0)-s_j(0) \| \} \ ,
\]
and choose $n_*$ such that 
\[
\| s_n(T)-s_{n_*}(0) \| = \min_{1 \leq k\leq N} \{ \| s_n(T)-s_k(0) \| \} \ .
\]
Then any $s_m(0)$ ($1 \leq m\leq N$) inside the ball $B_{\rho_nr_n}(s_{n_*}(0))$ implies 
that $s_m$ is admissible to follow $s_n$ in the segment sequence (cf. Definition \ref{segn}). 
One can also calculate interesting quantities as described at the end of subsection \ref{mard}.

The challenge in the numerical implementation of the segment description lies at the 
dimension of the absorbing set $\A$. Along each dimension, if we collocate $10$ points, and if the 
dimension of $\A$ is $n$, then we have total $10^n$ collocation points. So when the dimension 
$n$ is greater than $9$, the computation will be very challenging. The hope is that future 
computer technology may overcome this difficulty.

For low dimensional system like Lorenz equations, the numerical implementation of 
the segment description should not be difficult. For 2D Navier-Stokes equations with 
periodic boundary condition, numerical implementation of the segment description
for a Galerkin truncation of less than $6$ modes should not be difficult either. 

\section{Conclusion}

In this article, we introduced a segment description of turbulent solutions. The goal 
we are aiming at is to develop an effective description of turbulent solutions rather 
than the Reynolds average. Via the segment description, every turbulent solution can be 
approximated uniformly over an infinite time interval by a pseudo-solution built by segments of 
solutions. In dynamical systems, shadowing lemma addresses the converse of the above claim. 
In certain settings, shadowing lemma holds \cite{Li04} and for every pseudo-solution, one 
can find a true solution nearby. Even though the total number of segments in the 
segment description is finite, the number can be very large when the dimension of the phase 
space is large. This poses a challenge in numerical implementation of the segment description.

\end{document}